\documentclass[fleqn]{mat01}
\usepackage{times,mathtimy,amssymb,latexsym}
\begin{document}

\setcounter{page}{477} \firstpage{477}

\renewcommand\theequation{\thesection\arabic{equation}}

\newtheorem{theore}{Theorem}
\renewcommand\thetheore{\arabic{section}.\arabic{theore}}
\newtheorem{definit}[theore]{\rm DEFINITION}
\newtheorem{theor}{\bf Theorem}
\newtheorem{propo}[theore]{\rm PROPOSITION}
\newtheorem{lem}[theore]{\it Lemma}
\newtheorem{rem}[theore]{\it Remark}

\newtheorem{coro}[theore]{\rm COROLLARY}
\newtheorem{probl}[theore]{\it Problem}
\newtheorem{exampl}[theore]{\it Example}
\newtheorem{pot}[theore]{\it Proof of Theorem}

\def\nota{\trivlist \item[\hskip \labelsep{\it Notations.}]}

\title{Representations of homogeneous quantum L\'evy fields}

\markboth{V P Belavkin and L Gregory}{Representations of
homogeneous quantum L\'evy fields}

\author{V P BELAVKIN and L GREGORY}

\address{School of Mathematical Sciences, University of Nottingham,\\
\noindent Nottingham NG7 2RD, UK\\
\noindent E-mail: Viacheslav.Belavkin@nottingham.ac.uk;
dr.lee.gregory@gmail.com\\[1.2pc]
\noindent {\it Dedicated to K~B~Sinha on the \vspace{-1pc}event of
his 60th birthday}}

\volume{116}

\mon{November}

\parts{4}

\pubyear{2006}

\Date{}

\begin{abstract}
We study homogeneous quantum L\'{e}vy processes and fields with
independent additive increments over a noncommutative *-monoid.
These are described by infinitely divisible generating state
functionals, invariant with respect to an endomorphic injective
action of a symmetry semigroup. A~strongly covariant GNS
representation for the conditionally positive logarithmic
functionals of these states is constructed in the complex
Minkowski space in terms of canonical quadruples and isometric
representations on the underlying pre-Hilbert field space. This is
of much use in constructing quantum stochastic representations of
homogeneous quantum L\'{e}vy fields on It\^{o} monoids, which is a
natural algebraic way of defining dimension free, covariant
quantum stochastic integration over a space-time indexing set.
\end{abstract}


\maketitle

\section{Introduction}

The reconstruction theorem due to Kolmogorov is a celebrated
result in the theory of stochastic processes. It allows one to
build versions of a stochastic process $Z_{t}\hbox{\rm :}\
\Omega\mapsto\mathbb{R}^{d}$ in the narrow sense from a consistent
family of joint probabilities, called a stochastic process in the
wide sense, parametrised by $t$. There is a quantum multiparameter
generalisation of this theorem \cite{Belavkin85} allowing the
construction of quantum stochastic processes and fields given a
projective homogeneous system of multikernel maps describing the
process in the wide sense. The reconstruction of stationary
quantum weak Markov processes \cite{Bhat94} and quantum L\'{e}vy
fields \cite{Belavkin92c} parametrised by
$x\in\mathbb{R}_{+}\times\mathbb{R}^{d}$ instead of
$t\in\mathbb{R}_{+}$ are particular interesting cases of this
construction.

Our concern in this paper will be the quantum stochastic processes
and fields having independent increments in a weak differential
sense, called wide quantum L\'{e}vy fields. No assumption of
stationary increments will be made, instead we will consider
homogeneous processes with respect to a given action of a symmetry
semigroup in the field parameter space $X$. Classical
nonstationary L\'{e}vy processes are characterised by infinitely
divisible probability distributions, giving convolution hemigroups
\cite{Heyer95} of probability measures parametrised by the
intervals $\Delta= [t_{0},t_{1})$ of $\mathbb{R}_{+}$. In the
quantum setting, a commutative convolution hemigroup of states can
be defined by restricting the construction of \cite{Belavkin92c}
to one dimensional fields (that is, two parameter processes
indexed by intervals $\Delta= [t_{0},t_{1})$). There the coalgebra
structure determining the convolution of states was implicitly
defined as the commutative pointwise multiplication of generating
state functionals on the noncommutative monoid of a unital
$\star$-semigroup $\mathfrak{b}$, which will also be used here.
However, by enriching the enveloping semigroup algebra
$\mathbb{C}\mathfrak{b}$ of the $\star$-monoid $\mathfrak{b}$ with
a noncommutative coproduct one could use our results for general
L\'{e}vy processes and homogeneous fields as in the stationary,
one-parameter case \cite{Schurmann93}.

One can encode infinite divisibility without defining convolution
explicitly by using the necessary and sufficient exponential form
of the characteristic functional, called the generating state
functional in the quantum case. This generating state functional
is determined by its exponent, called the cumulant generating
state functional. This can be any Hermitian, absolutely
continuous, conditionally positive state functional on
$\mathfrak{b}$ which has unit $u\in\mathfrak{b}$ in its kernel.
The problem of reconstructing a quantum L\'{e}vy process in the
narrow sense is then reduced to finding a representation of the
cumulant generating state and exponentiating it in some sense.
This exponentiation is defined, similarly to the noncommutative
convolution, by a variety of independences (Boolean, tensor,
monotone and so on) which all coincide with the usual exponential
in the weak sense for the commutative coalgebra structure implicit
here. For non-commutative convolutions, this will result in the
noncommutative exponentiation depending on the choice of
independence, as in \cite{Franz05a}. Regardless, in all cases the
representation of the exponent is the\break same.

Here we present the first step by constructing `differential' type
representations associated with conditionally positive functionals
$\mathfrak{b}\rightarrow\mathbb{C}$ depending covariantly on the
field parameter $x\in X$ with respect to a semigroup of symmetries
acting both on $X$ and $\mathfrak{b}$. Later work will consider
the exponentiation mentioned above, giving stochastic integral
representations of covariant infinitely divisible positive
definite functionals. Covariant quantum dynamical semigroups are
an example of the noncommutative extension of this, and their
representations have been studied in \cite{Chakraborty01}. Our
main interest is when $\mathfrak{b}$ is obtained by a unitization
of a noncommutative It\^{o} algebra $\mathfrak{a}$ as a
parameterising algebra for the quantum stochastic differentials of
a quantum L\'{e}vy process as operator-valued processes with
independent increments in a quite general noncommutative\break
sense.

\section{Representations of homogeneous conditionally positive
functionals\\ on $\star$-semigroups}

Let $(X,\mathfrak{F},\mu)$ be a measurable space $X$ with a
$\sigma$-algebra $\mathfrak{F}$ and a positive $\sigma$-finite
atomless measure $\mu \hbox{\rm :}\
\mathfrak{F}\ni\Delta\mapsto\mu_{\Delta},\,\mu_{\mathrm{d}x}\equiv
\mathrm{d}x:=\mathrm{d}\mu(x)$, and let $\mathfrak{b}$ be a
semigroup with involution
\begin{equation*}
b\mapsto b^{\star},\quad (a\cdot c)^{\star}=c^{\star}\cdot
a^{\star},
\end{equation*}
and neutral element (unit) $u=u^{\star}$, $u\cdot b=b=b\cdot u$
for any $b\in\mathfrak{b}$. Typically $\mathfrak{b}$ will be a
unitization of a noncommutative It\^{o} $\star$-algebra
$\mathfrak{a}$, in which case
\begin{equation*}
a\cdot c=a+c+ac.
\end{equation*}
if $u$ is identified with zero, or simply write $a\cdot c=ac$ if
$\mathfrak{a}$ is realised as a $\star$-subalgebra of a unital
algebra by taking $u$ $=1 $. However in what follows one can take
any group with $u=1$ and $b^{\star}=b^{-1}$ or any
$\star$-submonoid of an operator algebra $\mathcal{B}$, a unit
ball of a unital $C^{*}$-algebra say, or even a filter (i.e. a
submonoid) of an idempotent, Boolean say, algebra $\mathfrak{B}$
with trivial involution $b^{\star}=b$.

Denote by $\mathfrak{m}$ the monoid of integrable step-maps
$g\hbox{\rm :}\ X\rightarrow \mathfrak{b}$, that is,
$\mathfrak{b}$-valued functions $x\mapsto g(x)$ having countable
images $g(X)=\{g(x)\hbox{\rm :}\ x\in X\}\subseteq\mathfrak{b}$
and integrable preimage $\Delta(b)=\{x\in X\hbox{\rm :}\
g(x)=b\}\in\mathfrak{F}$ in the sense$\,\mu _{\Delta(b)}<\infty$
for all $b\in\mathfrak{b}$ except $b=u$. We define on
$\mathfrak{m}$ an inductive structure of a $\star$-monoid with
pointwise defined operations $g^{\star}(x)=g(x)^{\star}$, $(f\cdot
h)(x)=f(x)\cdot h(x)$ and unit $e(x)=u$ for all $x\in X$,
\hbox{considering} $\mathfrak{m}$ as the union
$\cup\mathfrak{m}_{\Delta}$ of subsemigroups
$\mathfrak{m}_{\Delta}$ of functions $g\in\mathfrak{m}$ having
integrable\break supports
\begin{equation*}
\mathrm{supp}\,g=\{x\in X\hbox{\rm :}\ g(x)\neq u\}
\end{equation*}
in a $\Delta\in\mathfrak{F}$ with$\,\mu_{\Delta}<\infty$.

It is convenient to describe the $\star$-monoid $\mathfrak{b}$ by means
of a single Hermitian operation $a\star c=a\cdot c^{\star}$, satisfying
the relations
\begin{equation*}
b\star u=b,\quad u\star(u\star b)=b\quad\forall b\in\mathfrak{b}
\end{equation*}
defining $u=u^{\star}$ as right unit for the composition $\star$,
$b^{\star}$ as $u\star b$, and
\begin{equation*}
u\star((c\star b)\star a)=a\star((u\star b)\star c)
\end{equation*}
corresponding to $(a\cdot c^{\star})^{\star}=(a\star
c)^{\star}=c\star a=c\cdot a^{\star}$ and associativity of the
semigroup operation $a\cdot c$. This allows us to define both the
product and involution in a $\star$-monoid $\mathfrak{m}$ by a
single Hermitian binary operation $f\star h=g$, $\,g(x)=f(x)\star
h(x)$ with left unit $e\in\mathfrak{m}$ which recovers the
involution by $g^{\star}(x)=e\star g$ and the associative product
by $f\cdot h=f\star(e\star h)$ for all $f,h\in\mathfrak{m}$.

We introduce a semigroup $S$, called the \emph{symmetry
semigroup}, which has a measurable action on $X$ given by
injections $x\mapsto sx$, and denote for each $x\in sX$ its
preimage $s^{-1}x$ which is the unique element $x^{s}\in X$ such
that $sx^{s}=x$. The measure $\mu$ is assumed invariant under this
action in the sense that $\mu_{\Delta^{s}}=\mu_{\Delta}$ for each
$s\in S$ and any measurable $\Delta\subseteq sX$ with
$\Delta^{s}=s^{-1}\Delta$. We admit also an action $b\mapsto
b_{s}$ of $S$ on $\mathfrak{b}$ determined by a representation of
$S$ in the semigroup of unital injective $\star $-endomorphisms
$\theta_{s}\hbox{\rm :}\ \mathfrak{b}\rightarrow\mathfrak{b}$, so
that $\theta_{s}(u) =u$, $\theta_{s}(a\star b)=\theta_{s}(a)
\star\theta_{s}(b)$, and therefore
\begin{equation*}
(a\cdot c)^{s}=a^{s}\cdot c^{s},\quad b^{s\star}=b^{\star
s},\;u^{s}=u\;\;\;\forall a,b,c\in\theta_{s}(\mathfrak{b})
\end{equation*}
with respect to the inverse action $b\mapsto\theta_{s}^{-1}(b)
\equiv b^{s}$ of each $\theta_{s}$ on the $\star$-submonoid
$\mathfrak{b}_{s}=\theta_{s}(\mathfrak{b})$. These actions induce
a representation of $S$ on $\mathfrak{m}$ by the injective
$\star$-endomorphisms $g\mapsto g_{s}$,
\begin{equation*}
g_{s}(x)= \left\{
\begin{array}[c]{c@{\quad}c}
g (x^{s})_{s}, &x\in sX\\[.3pc]
u, & x\notin sX
\end{array}, \right.
\end{equation*}
obviously having the property that for any $s\in S$ and
$f,g\in\mathfrak{m}$ there exist unique
$f_{s},g_{s}\in\mathfrak{m}_{s}$ such that
\begin{equation*}
(f_{s}\star g_{s})^{s}(x):=(f_{s}(sx) \star g_{s}(sx))^{s} =
(f\star g)(x)
\end{equation*}
for all $x\in X$. Here we denoted by $\mathfrak{m}_{s}$ the
$\star$-submonoid of $\mathfrak{m}$ consisting of functions $g$
such that $g(x)\in\mathfrak{b}_{s}$ if $x\in sX$ and $g(x)=u$ if
$x\notin sX$, of which $g_{s}$ is a member.

We say that a complex functional $\varphi$ on $\mathfrak{m}$ is a
\emph{generating state functional} over the monoid $\mathfrak{m}$
(briefly, a \emph{state} over $\mathfrak{m}$), if the mapping
$\varphi\hbox{\rm :}\ \mathfrak{m}\rightarrow\mathbb{C}$ satisfies
the normalisation condition $\varphi(e)=1$ and positive
definiteness
\begin{equation}
\sum_{f,h\in\mathfrak{m}}\kappa_{f}\varphi(f\star h)
\kappa_{h}^{\ast} \geq0, \quad \forall\kappa_{g}\in
\mathbb{C}\hbox{\rm :}\ |\mathrm{supp}\,\kappa |<\infty,\label{one
a}
\end{equation}
where $|\cdot|$ denotes the cardinality of the set
$\mathrm{supp}\,\kappa =\{g\in\mathfrak{m}\hbox{\rm :}\
\kappa_{g}\neq0\}$. Every such function is lifted to a positive
normalised linear functional on the semigroup enveloping algebra
$\mathfrak{B}=\mathbb{C}\mathfrak{b}$. The state $\varphi$ is
called \emph{$S$-homogeneous} if $\varphi(g_{s})=\varphi(g)$ for
all $s\in S$ and $g\in\mathfrak{m}$.

Following \cite{Belavkin92c} we introduce on $\mathfrak{m}$ a
commutative and associative partial operation $f\sqcup\,h:=f\cdot
h$ for any functions $f,h\in\mathfrak{m}$ with disjoint supports
$\mathrm{supp}\,f\cap \mathrm{supp}\,h=\emptyset$. Thus the
defined map $\mathfrak{m}_{\Delta}
\times\mathfrak{m}_{\Delta^{\prime}}\rightarrow\mathfrak{m}_{\Delta}
\sqcup\mathfrak{m}_{\Delta^{\prime}}$ for any measurable disjoint
$\Delta,\Delta^{\prime}\in\mathfrak{F}$ is obviously lifted to the
tensor product $\mathbb{C}\mathfrak{m}_{\Delta}
\otimes\mathbb{C}\mathfrak{m}_{\Delta^{\prime}}$ of the enveloping
semigroup algebras of the $\star $-monoids $\mathfrak{m}_{\Delta}$
and $\mathfrak{m}_{\Delta^{\prime}}$. The operation $\sqcup$ is
well defined even for an infinite countable family $\{g_{n}\},
g_{n} \in\mathfrak{m}$ with mutually disjoint supports
$\Delta_{n}=\mathrm{supp}\,g_{n}$, by $\sqcup g_{n}(x)=g_{m}(x)$
for all $x\in\mathrm{supp}\,g_{m}$ and any $m$, otherwise $\sqcup
g_{n}(x)=u $ if $x\notin\sum\Delta_{n}$. Taking any $g
\in\mathfrak{m}$ and the partition
$\mathrm{supp}\,g=\sum\Delta_{n}$ into the co-images
$\Delta_{n}=\Delta(b_{n})$ where $b_{n} = g(x)$ for any $x
\in\Delta_{n}$, we see that any function $g\in\mathfrak{m}$ can be
written as $\sqcup g_{n}$, where $g_{n} = (b_{n})_{\Delta_{n}}$,
the $b_{n}$-valued indicator on $\Delta_{n}$. The $b$-valued
indicator of the subset $\Delta\subseteq X$ is defined in the
usual way: $b_{\Delta}(x)=b$ for all $x\in\Delta$ and
$b_{\Delta}(x)=u$ for $x\notin\Delta$.

We call a state $\varphi$ over $\mathfrak{m}$ \emph{chaotic} if it
satisfies the $\sigma$-multiplicativity condition
\begin{equation*}
\varphi\left(\bigsqcup_{n=1}^{\infty}g_{n}\right)
=\prod_{n=1}^{\infty}\varphi(g_{n}),
\end{equation*}
where $\prod_{n=1}^{\infty}\varphi(g_{n})=
\lim_{N\rightarrow\infty}\prod _{n=1}^{N}\varphi(g_{n})$ for any
functions $g_{n}\in\mathfrak{m}$ with pairwise disjoint supports:
$\mathrm{supp}\,g_{n}\cap\mathrm{supp} \,g_{m}=\emptyset$ for all
$n\neq m$. This condition is obviously fulfilled for $\varphi$ of
the exponential form $\varphi(g)=\hbox{e}^{\lambda(g)}$ with
\begin{equation}
\lambda(g)=\int l(x,g)\mathrm{d}x,\quad l(x,g)=l_{x}
(g(x)),\label{one b}
\end{equation}
which corresponds to absolute continuity (for all
$\Delta\in\mathfrak{F}$ we have
$\mu_{\Delta}=0\Rightarrow\lambda_{\Delta}(b)=0$) of the
$\sigma$-additive measure $\lambda_{\Delta}(b):=\lambda(
b_{\Delta})$ for each $b\in\mathfrak{b}$.

The family $\varphi_{\Delta}\hbox{\rm :}\
\mathfrak{b}\mapsto\mathbb{C}$ defined by any chaotic state
$\varphi$ as $\varphi_{\Delta}(b)=\varphi(b_{\Delta})$ is called
\emph{infinitely divisible} in the sense of the equality
$\varphi_{\Delta }(b)=\prod\varphi_{\Delta_{l}}(b)$ which also
holds in the limit of any integral sum sequence given by the
decomposition $\Delta=\Sigma\Delta_{i}$,
$\mu_{\Delta_{i}}\!\searrow\!0$ since
$\varphi_{\Delta_{i}}(b)\rightarrow1$ for any $b\in\mathfrak{b}$.
The function $\varphi_{\Delta}\hbox{\rm :}\ \mathfrak{b}
\rightarrow\mathbb{C}$ given by
\begin{equation}
\varphi_{\Delta}(b)=\exp\left\{
\int_{\Delta}l_{x}(b)\mathrm{d}x\right\} ,\label{one c}
\end{equation}
is clearly infinitely divisible in this sense.

Note that if the Radon--Nikodym derivative $l_{x}( b)
=\mathrm{d}\lambda( b) /\mathrm{d}x$ of the absolutely continuous
measure $\mathrm{d}\lambda( b) =\lambda_{\mathrm{d}x}( b) $ does
not depend on $x$, the states $\varphi_{\Delta}
(b)=\hbox{e}^{l(b)\mu_{\Delta}}\equiv\varphi^{\mu_{\Delta}}( b) $
form a continuous Abelian semigroup
\begin{equation*}
\{\varphi^{t}\hbox{\rm :}\
t\in\mathbb{R}^{+}\},\quad\varphi^{0}(b)=1,\quad[
\varphi^{r}\cdot\varphi^{s}]  (b)=\varphi^{r+s}(b)
\end{equation*}
with respect to the pointwise multiplication of $\varphi^{t}$.
\pagebreak

The $\lambda_{\Delta}$ introduced above is called the
\emph{cumulant generating state functional} and in general these
are defined as any function $b \mapsto\lambda_{\Delta}(b)$ which
is conditionally positive definite
\begin{equation}
\sum_{a,c\in\mathfrak{b}}\kappa_{a}\lambda_{\Delta}(a\star
c)\kappa_{c}^{\ast }\geq0,\,\forall\kappa\hbox{\rm :}\
|\mathrm{supp}\,\kappa|<\infty,\qquad\sum_{b\in
\mathfrak{b}}\kappa_{b}=0,\label{one d}
\end{equation}
such that $\lambda_{\Delta}(u)=0$, $\lambda_{\Delta}(b^{\star}
)=\lambda_{\Delta}(b)^{\ast}$ for any $b\in\mathfrak{b}$. They are
called \emph{$S$-homogeneous} if $\lambda(g_{s})=\lambda(g)$ for
all $g\in \mathfrak{m}$, which implies the $S$-homogeneity of the
generating state functional $\varphi$.

The following theorem shows that, along with some
differentiability conditions, the properties of a cumulant
generating state functional are necessary and sufficient
conditions for $\hbox{e}^{\lambda_{\Delta}}$ to be an
$S$-homogeneous infinitely divisible state. They are also
necessary and sufficient to allow the construction of a covariant
Minkowski space dilation. We assume that $X$ admits a net of
decompositions of the Vitali system in which
$\mu_{\Delta}\!\searrow\!0$, $\,x\in\Delta$, as
$\Delta\searrow\{x\}$.

\begin{theor}[\!]\label{T1} Consider an arbitrary functional
$\varphi_{\Delta}\hbox{\rm :}\ \mathfrak{b} \mapsto\mathbb{C},$
defined for any set $\Delta\in\mathfrak{F}$ of finite measure
$\mu_{\Delta}<\infty$ as $\varphi( b_{\Delta})$ by an
$S$-homogeneous functional $\varphi\hbox{\rm :}\
\mathfrak{m}\rightarrow\mathbb{C}$. The following are
equivalent{\rm :}
\begin{enumerate}
\renewcommand\labelenumi{\rm (\roman{enumi})}
\leftskip .4pc
\item $\varphi_{\Delta}$ is an infinitely divisible state
over $\mathfrak{b},$ and is an absolutely continuous
multiplicative measure in the sense that
$\mu_{\Delta}=0\Rightarrow\varphi_{\Delta}(b)=1$ for all
measurable $\Delta\in\mathfrak{F},\,b\in\mathfrak{b},$ and the
limit
\begin{equation}
\hskip -1.25pc
l_{x}(b)=\lim_{\Delta\downarrow\{x\}}\frac{1}{\mu_{\Delta}}
(\varphi_{\Delta }(b)-1)\label{one c'}
\end{equation}
exists in the Lebesgue--Vitali sense.

\item The functional $\lambda(g) = \ln\varphi(g)$ is defined{\rm ,} is
absolutely continuous in the sense that $\mu_{\Delta} = 0
\Rightarrow \lambda_{\Delta}(b) = 0$ for all measurable
$\Delta\in\mathfrak{F},\,b\in\mathfrak{b}${\rm ,} and is an
$S$-homogeneous cumulant generating functional with $S$-covariant
Radon--Nikodym derivative
\begin{equation*}
\hskip -1.25pc
l_{x}(b)=\lim_{\Delta\downarrow\{x\}}\frac{\lambda_{\Delta}( b)
}{\mu_{\Delta}}=l_{sx}(  b_{s})
\end{equation*}
for each $b\in\mathfrak{b}$.

\item There exists{\rm :}
\renewcommand\labelenumii{\rm (\arabic{enumii})}
\leftskip .5pc \vspace{6pt}

(1) an $S$-homogeneous integral $\star $-functional $\lambda(g)
=\int l(x,g)\,\mathrm{d}x = \ln \varphi(g)$ which has complex
$S$-covariant density $l( x,g) =l_{x}(g(x)) =l(sx,g_{s})$ such
that $l(x,g)^{\ast}\!=l(x,g^{\star})$ and whose values
$l(x,g)\!=0$\,\, for all $g(x)=u$ and $l(x,b_{\Delta})=l_{x}(b)$
with $x\in\Delta$ are independent of $\Delta;$

(2) a vector map $k\hbox{\rm :}\
g\mapsto\int^{\oplus}k(x,g)\,\mathrm{d} x$ into a pre-Hilbert
subspace $K\subseteq\int^{\oplus}K_{x}\mathrm{d}x$ of square
integrable functions $x\mapsto k( x,g) = k_{x}(g(x)) \in K_{x} $
which are \hbox{$S$-covariant}
\begin{equation*}
\hskip -2.5pc k(sx,g_{s})=k_{sx}(g(x)_{s}) =V_{s} k(x,g)
\end{equation*}
in terms of an isometric representation $s\mapsto V_{s}$ of $S$
with respect to the scalar products $\langle k|k\rangle \equiv\int
k_{x}^{\ast} k_{x}\mathrm{d}x=\Vert V_{s} k\Vert ^{2}$. It has
values $k(x,b_{\Delta })=k_{x}(b)\in K_{x}$ independent of
$\Delta\ni x$ and $k(x,b_{\Delta})=0$ if $x\notin\Delta$ such that
$k(x,g)=0$ if $g(x)=u$. The map $k$, together with the adjoint
functions $k^{\star}(x,g) =k( x,g^{\star})^{\ast}$ as the linear
functionals $k^{\star}(g) =\int^{\oplus}k^{\star} (x,g)
\mathrm{d}x\in K^{\ast}$, satisfies the condition
\begin{equation}
\hskip -2.5pc k^{\star}(f) k(h)=\lambda(f\cdot h) -\lambda(f)
-\lambda(h),\quad\forall f,h\in\mathfrak{m};\label{one e}
\end{equation}

(3) a unital $\ast$-representation $j\hbox{\rm :}\ g\mapsto
G:=\int^{\oplus}j(x,g)\mathrm{d} x,$
\begin{equation*}
\hskip -2.5pc j(g\cdot h)  =j(g)j(h),\;j(g^{\star})
=j(g)^{\ast},\;j(e) =I,
\end{equation*}
of the $\star$-monoid $\mathfrak{m}$ in the $\ast$-algebra of
decomposable operators $G\hbox{\rm :}\ K\ni
k\mapsto\int^{\oplus}j(x,g)k (x)\mathrm{d}x$ with
$j(x,b_{\Delta})=j_{x}(b)$ independent of $\Delta\ni x $ and
$j(x,b_{\Delta})=I_{x}$ if $x\notin\Delta${\rm ,} which are
$S$-covariant in the sense $V_{s}j(g) =j(g_{s}) V_{s}${\rm ,}
where $j( x,g_{s}) =j_{x}(g(x^{s})_{s})${\rm ,} $x\in sX${\rm ,}
otherwise $j( x,g_{s}) =I_{x}${\rm ,} satisfy the cocycle property
\begin{equation}
\hskip -2.5pc j(g)k(h)=k(g\cdot h) -k(g) ,\;\;\,k^{\star}(f)
j(g)=k^{\star}(f\cdot g) -k^{\star}(g) \label{one f}
\end{equation}
for all $f,g,h\in\mathfrak{m}$, and are continuous in $K$ with
respect to the polynorm
\begin{equation}
\hskip -2.5pc \Vert k\Vert^{h}= \left( \int\Vert
j(x,h)k(x)\Vert_{x} ^{2}\mathrm{d}x \right)^{1/2},\;\ \
\;h\in\mathfrak{m}.\label{one g}
\end{equation}
\item For each integrable $\Delta\in\mathfrak{F}$ there exists a
unital $\dagger$-representation $\mathbf{j}_{\Delta}\hbox{\rm :}\
\mathfrak{b}\rightarrow \mathfrak{b}( K_{\Delta}),$
\begin{equation*}
\hskip -1.25pc \mathbf{j}_{\Delta}(a\cdot b)=
\mathbf{j}_{\Delta}(a)
\mathbf{j}_{\Delta}(b),\;\;\mathbf{j}_{\Delta}
(b^{\star})=\mathbf{j}_{\Delta}(b)^{\dagger},
\mathbf{\;j}_{\Delta}(u)=\mathbf{I}_{\Delta}
\end{equation*}
in the algebra $\mathfrak{b}(K_{\Delta})$ of linear triangular
operators $\mathbf{L}=\mathbf{j}_{\Delta}(b)$. The integrals
\begin{equation*}
\hskip -1.25pc \lambda_{\Delta}(b) =
\int_{\Delta}l_{x}(b)\mathrm{d}x = \ln\varphi_{\Delta}(b), \qquad
j_{\Delta}(b) = \int^{\oplus}_{\Delta}j_{x}(b)\mathrm{d}x,
\end{equation*}
which correspond to $\sigma$-additivity and absolute continuity
with respect to the \hbox{$S$-homogeneous} measure $\mu$, appear
in the explicit form of $\mathbf{j}_{\Delta}$ given by
\begin{equation*}
\hskip -1.25pc \mathbf{j}_{\Delta}(b)= \begin{bmatrix} 1 &
k_{\Delta}^{\star}(b) & \lambda_{\Delta}(b)\\[.3pc] 0 & j_{\Delta}(b) &
k_{\Delta}(b)\\[.3pc] 0 & 0 & 1 \end{bmatrix} ,\quad
\mathbf{j}_{\Delta}(b)^{\dagger}= \begin{bmatrix} 1 &
k_{\Delta}(b)^{\ast} & \lambda_{\Delta}(b^{\star})\\[.3pc] 0 &
j_{\Delta}(b)^{\ast} & k_{\Delta}(b^{\star})\\[.3pc] 0 & 0 & 1
\end{bmatrix}
\end{equation*}
on the pseudo-Hilbert space $\mathbb{K}_{\Delta}=\mathbb{C}\oplus
K_{\Delta}\oplus\mathbb{C}$ defined by a pre-Hilbert space
$K_{\Delta}\subseteq \int_{\Delta}^{\oplus} K_{x}\mathrm {d}x$
with respect to the Minkowski scalar product
\begin{equation}
(\mathbf{k}|\mathbf{k}) _{\Delta}:=k_{-}^{\ast}k_{+}
+\int_{\Delta}\langle k_{x}|k_{x}\rangle
+k_{+}^{\ast}k_{-}\equiv\langle \mathbf{k}_{\Delta}^{\dagger}
,\mathbf{k}_{\Delta}\rangle .\label{one i}
\end{equation}
$\mathbf{L}^{\dagger}$ is the pseudo-Hermitian adjoint
$(\mathbf{k} |\mathbf{L}^{\dagger}\mathbf{k}) =(\mathbf{Lk}|
\mathbf{k})${\rm ,} $\mathbf{k}\in\mathbb{K}_{\Delta}\;$ and
$\lambda_{\Delta}(b) =\langle
\mathbf{e}^{\dagger},\mathbf{j}_{\Delta}( b) \mathbf{e}\rangle$
with respect to the row-vector $\mathbf{e}^{\dagger}=( 1,0,0)
\in\mathbb{K}_{\Delta }^{\dagger}$ of zero pseudo-norm $\langle
\mathbf{e}^{\dagger},\mathbf{e}\rangle _{\Delta}=0$. The family of
representations $\{ \mathbf{j}_{\Delta}\}$ is $S$-covariant in the
sense that there exists a representation $\mathbf{V}_{s}$ of $S$
in the pseudo-isometric operators
$\mathbb{K}_{\Delta}\rightarrow\mathbb{K}_{s\Delta}$ such that
$\mathbf{V}_{s}\mathbf{j}_{\Delta}(b)=\mathbf{j}_{s\Delta}(b_{s}
)\mathbf{V}_{s}$ for any $b\in\mathfrak{b}$ and integrable
$\Delta\subseteq X$. The pseudo-isometry $\mathbf{V}_{s}$ is
block-diagonal with $[\mathbf{V} _{s}]^{-}_{-} =
[\mathbf{V}_{s}]^{+}_{+} = 1${\rm ,} $[\mathbf{V}_{s}]^{\circ
}_{\circ} = V_{s}$ for some direct integral pseudo-isometry $V_{s}
\hbox{\rm :}\ K_{\Delta} \mapsto K_{s\Delta}${\rm ,} with all
other components zero.
\end{enumerate}
\end{theor}

We will prove that (iv) $\Rightarrow$ (iii) $\Rightarrow$ (ii)
$\Rightarrow$ (i) $\Rightarrow$ (iv) following \cite{Belavkin92c}
(see also \cite{Belavkin03}), adding homogeneity and covariance.

\noindent (iv) $\Rightarrow$ (iii). First let us define the
density $\mathbf{j}_{x}(b) = \lim_{\Delta\downarrow\{x\}}
\frac{\mathbf{j}_{\Delta}(b)}{\mu_{\Delta}}$ with the limit
understood in the appropriate sense for each of the blocks. Then
$\mathbf{j}_{x}(b)$ is a unital $\dagger$-representation of
$\mathfrak{b}$ in $\mathfrak{b}(K_{x})$. Thus we can define
$\mathbf{j}(x,g) = \mathbf{j}_{x}(g(x))$ which has direct integral
$\mathbf{j}(g) = \int_{X}^{\oplus}\mathbf{j}(x,g)\mathrm{d}x$.
Furthermore, $S$-homogeneity of $\mu$ gives the covariance
\begin{align*}
\mathbf{V}_{s}\mathbf{j}_{x}(b) &=
\lim_{\Delta\downarrow\{x\}}\frac
{\mathbf{V}_{s}\mathbf{j}_{\Delta}(b)}{\mu_{\Delta}} =
\lim_{\Delta \downarrow\{x\}}\frac{\mathbf{j}_{s\Delta}
(b_{s})\mathbf{V}_{s}}{\mu_{\Delta}}\\[.5pc]
&=
\lim_{\Delta\downarrow\{sx\}}\frac{\mathbf{j}_{\Delta}(b_{s})\mathbf{V}_
{s} }{\mu_{\Delta}} = \mathbf{j}_{sx}(b_{s})\mathbf{V}_{s}.
\end{align*}
Each of the components $\lambda(g)$, $k(g)$ and $j(g)$ of
$\mathbf{j}(g)$ is the direct integral of the component densities
$l_{x}(g(x))=l(x,g)$, $k_{x}(g(x))=k(x,g)$, $j_{x} (g(x))=j(x,g)$
existing due to the $\sigma$-additivity and absolute continuity
assumed in (iv). We can now use the properties of $\mathbf{j}_{x}$
to verify those required by (iii).

Starting with $\lambda(g) = \int^{\oplus}_{X}l(x,g)\mathrm{d}x$ we
have $l(x,g)= \langle\mathbf{e}^{\dagger},\mathbf{j}(x,g)
\mathbf{e} \rangle$ so that $l(x,g^{\star})=
\langle\mathbf{e}^{\dagger}, \mathbf{j}(x,g^{\star}) \mathbf{e}
\rangle = \langle\mathbf{e}^{\dagger},\mathbf{j}(x,g)^{\dagger}
\mathbf{e} \rangle= l(x,g)^{\ast}$ and therefore $l(x,e) = \langle
\mathbf{e}^{\dagger},\mathbf{j}(x,u) \mathbf{e} \rangle= 0$. Hence
$l(x,g)=l_{x}(u)=0$ for any $x $ in the kernel of $g$, and also
$l(x,b_{\Delta })$ is independent of $\Delta\ni x$. The covariance
of $\mathbf{j}_{x}$ gives
\begin{equation*}
l(x,g) = \langle\mathbf{e}^{\dagger},\mathbf{j}_{x}(g) \mathbf{e}
\rangle= \langle\mathbf{e}^{\dagger}, \mathbf{V}_{s}
\mathbf{j}_{sx} (g(x)_{s})\mathbf{V}_{s} \mathbf{e} \rangle=
l(sx,g_{s})
\end{equation*}
as required. The $S$-homogeneity $\lambda(g)=\lambda(g_{s})$ follows
straight away by integrating and noting that $l(x,g_{s})=0$ for $x
\notin sX$.

The vector map density $k(x,g)$ is given by the density
$[\mathbf{j}_{x}(g(x))]^{\circ}_{+}$. Due to the unitality of
$\mathbf{j}_{x} $, $k(x,g(x)) = 0$ for $x$ in the kernel of $g$,
and clearly $k(x,b_{\Delta})$ is independent of $\Delta\ni x$. To
prove condition $(\ref{one e})$ consider
\begin{equation}
\lambda(f\cdot g)= \langle\mathbf{e}^{\dagger} , \mathbf{j}(f
\cdot g) \mathbf{e} \rangle= \lambda(g) + k^{\star}(f)k(g) +
\lambda(f).\label{one j}
\end{equation}
This also shows that $\Vert k \Vert^{2} = k^{\star}(f)k(g) <
\infty$ so that $k(x,g)$ is square integrable as required.
Covariance again follows trivially by taking the appropriate
component of the covariance for $\mathbf{j}_{x}$:
\begin{equation*}
V_{s}k(x,g) = [\mathbf{V}_{s}\mathbf{j} (x,g)]^{\circ}_{+} =
[\mathbf{j}(sx,g_{s}) \mathbf{V}_{s}]^{\circ}_{+} = k(sx,g_{s}).
\end{equation*}
Extending $V_{s}$ to $K$ by its direct integral $V_{s}k(g) =
\int^{\oplus}_{X}V_{s}k(x,g)\mathrm{d}x$ also gives a linear
pseudo-isometry.

That the decomposable linear operator
$j(g)=[\mathbf{j}(g)]^{\circ}_{\circ}$ on $K$ is a unital
$\ast$-representation of $\mathfrak{m}$ follows from the upper
triangular form of $\mathbf{j}(g)$ and its $\dagger
$-multiplicativity and unitality. The cocycle property $(\ref{one
f})$ is seen straightforwardly by calculating $k(g \cdot h) =
[\mathbf{j} (g)\mathbf{j}(h)]^{\circ}_{+}$. Completing $K$ with
respect to the seminorms $(\ref{one g})$ obviously makes $j$
continuous. Covariance of $j$ again follows simply by taking the
appropriate component of the covariance condition for
$\mathbf{j}$.

(iii) $\Rightarrow$ (ii). It is immediate that the absolutely
continuous measure $\lambda_{\Delta}(b)$ satisfies the conditions
$\lambda_{\Delta }(b^{\star})=\lambda_{\Delta}(b)^{\ast}$ and
$\lambda_{\Delta}(u)=0 $, since the functional $l_{x}$ satisfies
the respective conditions. Also the condition $l_{x}(u)=0$ along
with the integral form of $\lambda_{\Delta}$ ensure its
$\sigma$-additivity on disjoint integrable subsets of $X$, and so
is the cumulant generating functional of an infinitely divisible
state if it is conditionally positive definite. The
\hbox{conditional} positivity (\ref{one d}) follows from (\ref{one
e}) and the positive definiteness of the inner product\break on
$K$:
\begin{align*}
\sum_{f,g\in\mathfrak{m}}\kappa_{f}\lambda(f\star g)
\kappa_{g}^{\ast} &= \sum_{f,g\in\mathfrak{m}}\kappa_{f} (
\lambda(g^{\star}) + \lambda(f) +
k^{\star}(f)k(g^{\star}) )  \kappa_{g}^{\ast}\\[.5pc]
&= \sum_{f,g\in\mathfrak{m}}\kappa_{f} k^{\star}(f)k(g^{\star})
\kappa _{g}^{\ast} \geq0
\end{align*}
for $\sum_{f}\kappa_{f} = 0$. $S$-homogeneity of $\lambda(g)$ and
$l_{x}$ are trivial.

(ii) $\Rightarrow$ (i). If the function $\lambda_{\Delta}(b)$ is a
(complex) absolutely continuous measure, then
$\varphi_{\Delta}(b)=\exp\{\lambda _{\Delta}(b)\}$ has the
property $\varphi_{\sqcup\Delta_{l}}(b)=\prod
\varphi_{\Delta_{l}}(b)$ of infinite divisibility. Moreover the
limit (\ref{one c'}) exists, and by virtue of
$\varphi_{\Delta}(b)\rightarrow1$ as $\Delta\!\downarrow\!\{x\}$
it coincides with the Radon--Nikod\'{y}m derivative
$l_{x}(b)=\mathrm{d}\ln\varphi(b)/\mathrm{d}x$ as the limit of the
quotient $\lambda_{\Delta}(b)/\mu_{\Delta}$ over a net of subsets
$\Delta\ni x$ of the system of Vitali decompositions of the
measurable space $X$. $S$-homogeneity of $\varphi_{\Delta}(b)$ is
trivial. For any integrable $\Delta$ the function
$b\mapsto\varphi_{\Delta}(b)$ can be shown to be positive in the
sense of (\ref{one a}) by considering the conditioned complex
function $b\mapsto \kappa^{\circ}_{b}$ defined as
$\kappa^{\circ}_{u} = \kappa_{u} -\sum_{b \in\mathfrak{b}}
\kappa_{b}$, $\kappa^{\circ}_{b} = \kappa_{b}$ for all $b \neq u
\in\mathfrak{b}$ for an arbitrary finitely supported complex
function $b\mapsto\kappa_{b}$. Hence
$\sum_{b\in\mathfrak{b}}\kappa^{\circ}_{b}=0$ so by conditional
positivity
\begin{align*}
0\leq\sum_{a,c\in\mathfrak{b}}\kappa^{\circ}_{a}\lambda_{\Delta}(a\star
c)\kappa^{\circ\ast}_{c} =
\sum_{a,c\in\mathfrak{b}}\kappa_{a}(\lambda _{\Delta}(a\star
c)-\lambda_{\Delta}(a)-\lambda_{\Delta}(c^{\star}))\kappa
_{c}^{\ast},
\end{align*}
where we have taken into account the fact that
$\lambda_{\Delta}(u)=0$. Since the exponential of any
positive-definite kernel is a positive definite kernel, we have
for any $\Delta$,
\begin{align*}
&\sum_{a,c\in\mathfrak{b}}\kappa_{a}^{\ast}\exp\{\lambda_{\Delta}(a\star
c)\}\kappa_{c}\\[.5pc]
&\quad\, = \sum_{a,c\in\mathfrak{b}}\kappa_{\Delta}^{a\ast}\exp
\{\lambda_{\Delta}(a\star
c)-\lambda_{\Delta}(a)-\lambda_{\Delta}(c^{\star})\}
\kappa_{\Delta}^{c}\geq0,
\end{align*}
where $\kappa_{\Delta}^{b} = \kappa_{b}\exp\{\lambda_{\Delta}(b)\}$ and we
have used $\lambda_{\Delta}(b^{\star}) = \lambda_{\Delta}(b)^{\ast}$.

(i) $\Rightarrow$ (iv). Since $\varphi_{\Delta}$ is an infinitely
divisible state on $\mathfrak{b}$ and
$\varphi_{\Delta}(b)\rightarrow1$ for all $b$ as
$\mu_{\Delta}\rightarrow0$, the limit $l_{x}(b)$ is defined as the
logarithmic derivative
$\mu_{\mathrm{d}x}^{-1}\ln\varphi_{\mathrm{d}x}(b)$ of the measure
$\lambda_{\Delta}(b)=\ln\varphi_{\Delta}(b)$ in the Radon--Nikodym
sense. Consequently, the function $x\mapsto l_{x}(b)$ is
integrable and almost everywhere satisfies the conditions
$l_{x}(a\star c)^{\ast}=l_{x}(c\star a)$, $\,l_{x}(u)=0$ and
\begin{equation*}
\sum_{b\in\mathfrak{b}}\kappa_{b}=0\;\Rightarrow(\kappa^{\prime}|
\kappa)_{x}:=\sum_{a,c\in\mathfrak{b}}\kappa_{a}l_{x}(a\star
c)\kappa _{c}^{\ast}\geq0
\end{equation*}
for all $\kappa$ such that $|\mathrm{supp}\,\kappa|<\infty$, which
can easily be verified directly for the difference derivative
$l_{\Delta}(b)=(\varphi _{\Delta}(b)-1)/\mu_{\Delta}$ passing to
the limit $\Delta\!\!\downarrow\!\!\{x\}$. In addition
$\int_{\Delta}l_{x}(b)\mathrm{d}x=\ln\varphi_{\Delta}(b) =
\lambda_{\Delta}(b)$ by absolute continuity, and since
$(b_{\Delta})_{s}(x) = (b_{s\Delta}(x))_{s}$, $S$-homogeneity
becomes $\varphi_{s\Delta}(b) = \varphi_{\Delta}(b^{s})$ for
$\Delta\subset sX$, giving $l_{sx}(b)=l_{x} (b^{s})$ for $x \in
sX$ in the limit.

We consider the space $\mathfrak{B}$ of complex functions
$\kappa=( \kappa_{b}) _{b\in\mathfrak{b}}$ on $\mathfrak{b}$ with
finite supports $\{b\in\mathfrak{b}\hbox{\rm :}\
\kappa_{b}\neq0\}$ as a unital $\star$-algebra with respect to the
product $\kappa^{\prime}\cdot\kappa$ defined as
$\kappa^{\prime}\star\kappa^{\star}$ by the Hermitian convolution
\pagebreak
\begin{equation*}
(\kappa^{\prime}\star\kappa)_{b}=\sum_{a\star
c=b}\kappa_{a}^{\prime}
\kappa_{c}^{\ast},\;\;\;\;\;\,\delta_{u}\star\kappa=\kappa^{\star
},\,\;\;\;\;\kappa\star\delta_{u}=\kappa
\end{equation*}
with right identity $\delta_{u}$. Here $\delta_{a}=(\delta
_{a,b})_{b\in\mathfrak{b}}$ is the Kronecker delta and it defines
a \hbox{$\star$-representation} $a\mapsto\delta_{a}$ of the monoid
$\mathfrak{b}$ in $\mathfrak{B}$,
\begin{equation*}
\delta_{a}\star\delta_{c}=\delta_{a\star c},\quad\delta_{u}\star\delta
_{b}=\delta_{b},\quad\delta_{b}\star\delta_{u}=\delta_{b^{\star}},
\end{equation*}
with respect to the involution $\kappa^{\star}=( \kappa_{b^{\star}
}^{\ast})  _{b\in\mathfrak{b}}$. The linear subspace $\mathfrak{A}
\subset\mathfrak{B}$ of distributions $\kappa$ such that the sum
$\kappa _{-}:=\sum_{b\in\mathfrak{b}}\kappa_{b}$ equals zero, is a
$\star$-ideal since
\begin{equation*}
\sum_{b\in\mathfrak{b}}(\kappa^{\prime}\star\kappa)_{b}=\sum_{b\in
\mathfrak{b}}\sum_{a\star
c=b}\kappa_{a}^{\prime}\kappa_{c}^{\ast}=\sum
_{a\in\mathfrak{b}}\kappa_{a}^{\prime}\sum_{c\in\mathfrak{b}}\kappa_{c}^{\ast
}=0.
\end{equation*}

Let us equip $\mathfrak{B}$ for every $x\in X$ with the\ Hermitian
form $(\kappa^{\prime}|\kappa)_{x}$ of the kernel $l_{x}( a\star
c) $ which is positive on $\mathfrak{A}$ and can be written in
terms of the kernel $\langle \delta_{a},\delta_{c}^{\star}\rangle
_{x}^{\circ }=l_{x}(a\star c)-l_{x}(a)-l_{x}(c^{\star})$ as
\begin{equation*}
(\kappa^{\prime}|\kappa)_{x}=\kappa_{-}^{\prime}\kappa_{+}^{\ast
}+\langle \kappa^{\prime},\kappa^{\star}\rangle _{x}^{\circ}
+\kappa_{+}^{\prime}\kappa_{-}^{\ast},
\end{equation*}
where $\kappa^{x}_{+}:=\sum_{b}\kappa_{b}l_{x}( b) $. We notice
that the Hermitian form
\begin{equation*}
\langle \kappa^{\prime\star}|\kappa^{\star}\rangle _{x}^{\circ
}:=\sum_{a.c\in\mathfrak{b}}\kappa_{a}^{\prime}\langle \delta_{a}
,\delta_{c}\rangle _{x}^{\circ}\kappa_{c}^{\star}\equiv\langle
\kappa^{\prime},\kappa^{\star}\rangle _{x}^{\circ}
\end{equation*}
is non-negative if $\kappa_{-}=0$ or $\kappa_{-}^{\prime}=0$ as
$\langle
\kappa,\kappa^{\star}\rangle_{x}^{\circ}=\sum\kappa_{a}\langle
\delta_{a},\delta_{c}^{\star}\rangle_{x}^{\circ}\kappa_{c}^{\ast}\geq
0$, coinciding with $(\kappa^{\prime}|\kappa)_{x}$. Since
$(\kappa^{\prime}|\kappa)_{x}=\sum_{b}(\kappa^{\prime}\star\kappa)_{b}
l_{x}(b)$, the form $(\kappa^{\prime}|\kappa)_{x}$ has right
associativity property
\begin{equation*}
(\kappa^{\prime}\cdot\kappa|\kappa)_{x}=(\kappa^{\prime}|\kappa
\star\kappa)_{x}=(\kappa^{\prime}|\kappa\cdot\kappa^{\star})_{x},
\end{equation*}
for all $\kappa,\kappa^{\prime}\in\mathfrak{B}$, and therefore its
kernel $\mathfrak{R}_{x}=\{ \kappa\hbox{\rm :}\
(\kappa^{\prime}|\kappa)_{x} =0, \forall\kappa^{\prime}\} $ is the
right ideal
\begin{equation*}
\mathfrak{R}_{x}=\{\kappa^{\prime}\in\mathfrak{B}\hbox{\rm :}\
(\kappa^{\prime}\cdot
\kappa|\kappa)_{x}=0,\,\forall\kappa\in\mathfrak{B}\}
\end{equation*}
belonging to $\mathfrak{A}$. We factorise $\mathfrak{B}$ by this
right ideal putting $\kappa\approx0$ if
$\kappa\in\mathfrak{R}_{x}^{\star}:=\{\kappa^{\star}\hbox{\rm :}\
\kappa\in\mathfrak{R}_{x}\} $ and denoting the equivalence classes
of the left factor-space $\mathbb{K}_{x} = \mathfrak{B}
/\mathfrak{R}_{x}^{\star}$ as the ket-vectors
$|\kappa)=\{\kappa^{\prime }\hbox{\rm :}\
\kappa^{\prime}-\kappa^{\star}\in\mathfrak{R}_{x}^{\star}\}$. The
condition $\kappa\in\mathfrak{R}_{x}$ means in particular that
$\kappa_{x}^{-} :=(\delta_{u}|\kappa)_{x}=0$, and therefore
\begin{equation*}
(\kappa|\kappa)_{x}=\sum_{a,c\in\mathfrak{b}}\kappa_{a}\langle\delta
_{a},\delta_{c}^{\star}\rangle^{\circ}
_{x}\kappa_{c}^{\ast}=\langle \kappa^{\circ}|\kappa^{\circ}\rangle
_{x}=0,
\end{equation*}
where $\kappa^{\circ}=( \kappa_{b}^{\circ}) _{b\in\mathfrak{b}}$
denotes an element of $\mathfrak{A}$ obtained as
$\kappa_{b}^{\circ} =\kappa_{b}^{\star}$ for all $\,b\neq u$ and
$\,\kappa_{u}^{\circ}=\kappa
_{u}^{\star}-\sum_{b\in\mathfrak{b}}\kappa_{b}^{\star}$ such that
$\langle \kappa^{\circ}|\kappa^{\circ}\rangle _{x}=\langle
\kappa,\kappa^{\star}\rangle _{x}^{\circ}$. Therefore it follows
also that $\kappa^{+}:=\sum\kappa_{b}^{\star}$ is also zero for
any $\kappa \in\mathfrak{R}_{x}$ since
\begin{equation*}
0=(\kappa^{\prime}|\kappa)_{x}=\kappa_{-}^{\prime}\kappa_{+}^{x\ast
}+\langle \kappa^{\prime},\kappa^{\star}\rangle _{x}^{\circ}
+\kappa_{+}^{\prime
x}\kappa_{-}^{\ast}=\kappa_{-}^{\ast}=\kappa^{+}
\end{equation*}
for any $\kappa^{\prime}\in\mathfrak{B}$ with $\kappa_{+}^{\prime
x}=1$ by virtue of $\kappa_{+}^{x \ast}=\kappa^{-}=0$ and also due
to the Schwartz inequality $(\kappa^{\prime}|\kappa)_{x}=\langle
\kappa^{\prime },\kappa^{\star}\rangle^{\circ}_{x}=0$. This allows
us to represent the left equivalence classes $|\kappa)_{x}$ by the
columns $\mathbf{k}=[k^{\mu}]$ with $k^{\mp}=\kappa^{\mp}$ and
$k^{\circ}=|\kappa^{\circ }\rangle_{x}$ in the Euclidean component
$K_{x} \subset \mathbb{K}_{x}$ as the subspace of the left
equivalence classes $|\kappa ^{\circ}\rangle_{x}
=|\kappa_{\circ})_{x}$ of the elements $\kappa_{\circ }=(
\kappa_{b}-\delta_{u,b}\kappa_{-}) _{b\in\mathfrak{b}}
\in\mathfrak{A}$ such that
$\kappa_{\circ}^{\star}=\kappa^{\circ}$. These columns are
pseudo-adjoint to the rows ${\pmb k}=(k_{-},k_{\circ} ,k_{+})$ as
the right equivalence classes $_{x}(\kappa:=|\kappa)_{x}^{\dagger
}\in\mathfrak{B}/\mathfrak{R}_{x}$ with $\,k_{\pm}=\kappa_{\pm}$
and $k_{\circ}= \, _{x}(\kappa_{\circ}$ defining the indefinite
product in terms of the canonical pairing
\begin{equation*}
k_{\cdot}k^{\cdot}=k_{-}k^{-}+\langle k_{\circ},k^{\circ}\rangle
_{x} +k_{+}k^{+}=({\pmb k}|{\pmb k}^{\dagger}) _{x} ,
\end{equation*}
where $k^{\circ}=k_{\circ}^{\ast}\in K_{x}$,
$k^{\pm}=k_{\mp}^{\ast }\in\mathbb{C}$ with respect to the
Euclidean scalar product $\langle k_{\circ},k^{\circ}\rangle _{x}
=\langle k_{\circ}^{\ast}| k^{\circ}\rangle_{x}$ of the Euclidean
space $K_{x}=\{ k^{\circ}=|\kappa^{\circ}\rangle_{x} \hbox{\rm :}\
\kappa_{\circ}\in\mathfrak{A}\}$, and
\begin{equation*}
\kappa_{+}^{x\ast}=\sum_{b\in\mathfrak{b}}l_{x}(b^{\star})\kappa^{b}
=\kappa_{x}^{-},\quad \kappa_{-}^{\ast}=\sum_{b\in\mathfrak{b}}
\kappa_{b}^{\ast}=\kappa^{+}.
\end{equation*}

We notice that the representation $\delta_{\cdot}\hbox{\rm :}\
\mathfrak{b}\ni b\mapsto\delta_{b}$ is Hermitian:
\begin{equation*}
(\kappa\cdot\delta_{b}|\kappa)_{x}=\sum_{b\in\mathfrak{b}}l(b)(\kappa
\cdot\delta_{b}\star\kappa)_{b}=(\kappa|\kappa\cdot\delta_{b^{\star}})_{x},
\end{equation*}
and that it is well defined as right representation on
$\mathfrak{B} /\mathfrak{R}_{x}$ (or left representation on
$\mathfrak{B}/\mathfrak{R} _{x}^{\star}$) since
$\kappa\cdot\delta_{b}\in\mathfrak{R}_{x}$ if $\kappa
\in\mathfrak{R}_{x}$:
\begin{equation*}
(\kappa|\kappa)_{x}=0,\quad
\forall\kappa\in\mathfrak{B}\Rightarrow(\kappa
\cdot\delta_{b}|\kappa)_{x}=(\kappa|\kappa\star\delta_{b})_{x}
=0,\quad \forall\kappa\in\mathfrak{B}.
\end{equation*}

This allows us to define for each $b\in\mathfrak{b}$ an operator
$_{x}(\kappa\mathbf{j}_{x}(b)= \, _{x}(\kappa\cdot\delta_{b}$ such
that $\mathbf{j} _{x}(b^{\star})=\mathbf{j}_{x}(b)^{\dagger}$ with
the componentwise action
\begin{align*}
(\kappa\cdot\delta_{b})_{-} &  =\kappa_{-},\;\;(\kappa\cdot\delta_{b})_{\circ
}=\kappa_{-}(\delta_{b}-\delta_{u})+\kappa_{\circ}\cdot\delta_{b},\\[.3pc]
(\kappa\cdot\delta_{b})_{+}^{x} &
=\kappa_{-}l_{x}(b)+(\kappa_{\circ}
|\delta_{b^{\star}}-\delta_{u})_{x} +\kappa_{+}^{x},
\end{align*}
given as the right multiplications ${\pmb k}\mapsto{\pmb k}
\mathbf{L}$, ${\pmb k}\mapsto{\pmb k}\mathbf{L}^{\dagger}$ of the
triangular matrices $\mathbf{L} \equiv\mathbf{j}_{x}(b)$,
$\mathbf{L}^{\dagger} \equiv\mathbf{j}_{x}(b^{\star})$,
\begin{align*}
\mathbf{L}=
\begin{bmatrix}
1 & j_{\circ}^{-}(x,b) & j_{+}^{-}(x,b)\\[.2pc]
0 & j_{\circ}^{\circ}(x,b) & j_{+}^{\circ}(x,b)\\[.2pc]
0 & 0 & 1
\end{bmatrix}
,\quad \mathbf{L}^{\dagger} =\left[
\begin{array}[c]{ccc}
1 & j_{+}^{\circ}(x,b)^{\ast} & j_{+}^{-}(x,b)^{\ast}\\[.2pc]
0 & j_{\circ}^{\circ}(x,b)^{\ast} & j_{\circ}^{-}(x,b)^{\ast}\\[.2pc]
0 & 0 & 1
\end{array}
\right]
\end{align*}
by the rows ${\pmb k}=(k_{-},k_{\circ},
k_{+})\in\mathbb{K}^{\dagger}_{x}$ (or as the left multiplications
$\mathbf{Lk}$, $\mathbf{L}^{\dagger }\mathbf{k}$ by columns
$\mathbf{k}\in\mathbb{K}_{x}$). Here
\begin{align*}
j_{+}^{-}(x,b) &=
l_{x}(b),\,\;\;\;\;_{x}(\kappa_{\circ}j_{\circ}^{\circ
}(x,b)=\,_{x}(\kappa_{\circ}\cdot\delta_{b}= \,
_{x}(\kappa_{\circ} j_{x}( b),\\[.3pc]
\,\;\;j_{+}^{\circ}(x,b^{\star})  &=\delta_{b}^{\star}\rangle_{x}
=k_{x}(  b^{\star}) = k_{x}^{\star}(  b) ^{\ast}=\,
_{x}\langle\delta_{b}^{\star}|^{\ast}=j_{\circ}^{-}(x,b)^{\ast},
\end{align*}
where $\delta_{b}^{\star}\rangle_{x}=|\delta_{b}-\delta_{u})_{x}$
and $\,L_{-\nu}^{\star\mu}=L_{-\mu}^{\nu\ast}$ is pseudo-Euclidean
conjugation of the triangular matrix $\mathbf{L}=[L_{\nu}^{\mu}]$
corresponding to the map ${\pmb k}\mapsto{\pmb k}^{\dagger}$ into
the adjoint columns  $\mathbf{k}=[k^{\mu}]$ with the components
$k^{\mu}=k_{-\mu}^{\ast}$ given by the pseudo-metric tensor
$g^{\mu\nu}=\delta_{-\nu}^{\mu }=g_{\mu\nu}$:
\begin{align*}
\begin{bmatrix}  b_{-}^{-} & b_{\circ}^{-} & b_{+}^{-}\\[.3pc]  0 &
b_{\circ}^{\circ} & b_{+}^{\circ}\\[.3pc]  0 & 0 & b_{+}^{+}
\end{bmatrix}^{\star} &= \begin{bmatrix}  0 & 0 & 1\\[.3pc]  0 & I & 0\\[.3pc]
1 & 0 & 0  \end{bmatrix}  \begin{bmatrix}  b_{-}^{-} &
b_{\circ}^{-} & b_{+}^{-}\\[.3pc]  0 & b_{\circ}^{\circ} &
b_{+}^{\circ}\\[.3pc]  0 & 0 & b_{+}^{+} \end{bmatrix}  ^{\ast}
\begin{bmatrix}  0 & 0 & 1\\[.3pc]  0 & I & 0\\[.3pc]  1 & 0 & 0
\end{bmatrix}\\[.5pc]
&= \begin{bmatrix}  b_{+}^{+\ast} &
b_{+}^{\circ\ast} & b_{+}^{-\ast}\\[.3pc]  0 & b_{\circ}^{\circ\ast} &
b_{\circ}^{-\ast}\\[.3pc]  0 & 0 & b_{\circ}^{-\ast}  \end{bmatrix}.
\end{align*}
Constructing the direct integral space $\mathbb{K} = \int^{\oplus}
_{X}\mathbb{K}_{x}\mathrm{d}x$ allows the definition of
$\mathbf{j}_{\Delta  }(b)$ as the block-wise direct integral
$\int^{\oplus}_{\Delta}\mathbf{j} _{x}(b)\mathrm{d}x$. Thus we can
write the constructed canonical $\dagger$-representations
$\mathbf{j}_{\Delta}(b)=[j_{\nu}^{\mu}(\Delta, b)]$ of the  monoid
$\mathfrak{b}$ in the pseudo-Euclidean space $\mathbb{K}$ of
columns  $\mathbf{k}=[k^{\mu}]$ in terms of the usual matrix
multiplication
\begin{align*}
&\begin{bmatrix} 1 & k^{\star}_{\Delta}(a) & \lambda_{\Delta}(a)\\[.3pc]
0 & j_{\Delta}(a) & k_{\Delta}(a)\\[.3pc]  0 & 0 & 1  \end{bmatrix}
\begin{bmatrix}  1 & k_{\Delta}^{\star}(b) & \lambda_{\Delta}(b)\\[.3pc]
0 & j_{\Delta}(b) & k_{\Delta}(b)\\[.3pc]  0 & 0 & 1  \end{bmatrix}\\[.7pc]
&\quad\, =
\begin{bmatrix}  1 &
k_{\Delta}^{\star}(b)+k_{\Delta}^{\star}(a)j_{\Delta}(b) & \lambda
_{\Delta}(b) +
k_{\Delta}^{\star}(a)k_{\Delta}(b)+\lambda_{\Delta}(a)\\[.3pc]  0 &
j_{\Delta}(a)j_{\Delta}(b) &
j_{\Delta}(a)k_{\Delta}(b)+k_{\Delta}(a)\\[.3pc]  0 & 0 & 1
\end{bmatrix}.
\end{align*}
Clearly this realises the conditionally positive function
$\lambda_{\Delta  }(b)$ as the value of the vector form
\begin{equation*}
\mathbf{e}^{\dagger}\mathbf{j}_{\Delta}(b)\mathbf{e}=e_{\mu}j_{\nu}^{\mu
}(\Delta,b)e^{\nu}=j_{+}^{-}(\Delta,b)=\lambda_{\Delta}(b)
\end{equation*}
with the column $\mathbf{e}=[\delta_{+}^{\mu}] ={\bf e}^{\dagger}$
the adjoint to the row ${\bf e}=(  1,0,0)  $ of zero pseudonorm
$\mathbf{e}^{\dagger}\mathbf{e}=e_{\mu}e^{\mu}=0$ for each $x$.
A~representation $V_{s}$ of the symmetry semigroup $S$ is defined
on $\mathfrak{B}$ as
\begin{equation*}
(V_{s}\kappa)_{b} = \sum_{a \in\mathfrak{b}}
\kappa_{a}\delta_{a_{s},b} =  \biggl\{
\begin{array}{c@{\quad}l}  \kappa_{b^{s}}, & \text{if } \ b
\in\mathfrak{b}_{s}\\[.3pc]  0, & \text{if } \ b \notin\mathfrak{b}_{s}
\end{array}.
\end{equation*}
This map is $(x,sx)$-isometric on $\mathfrak{B}$ in the sense that
\begin{align*}
(V_{s}\kappa^{\prime}\vert V_{s}\kappa)_{sx} = \sum_{a,b
\in\mathfrak{b} }\kappa^{\prime}_{a}l_{sx}(a_{s} \star
b_{s})\kappa_{b}^{\ast} = \sum_{a,b
\in\mathfrak{b}}\kappa^{\prime}_{a}l_{x}(a \star
b)\kappa^{\ast}_{b} =  (\kappa^{\prime}\vert\kappa)_{x}.
\end{align*}
The pseudo-adjoint $V_{s}^{\star}$ is well defined as a surjection
from $\mathfrak{B}_{s}=\mathbb{C}\mathfrak{b}_{s}$ onto
$\mathfrak{B}$ so that  every $\kappa^{\prime}\in\mathfrak{B}$ can
be written as $V^{\star}_{s}\kappa^{\prime\prime}$ for some
$\kappa^{\prime\prime}\in\mathfrak{B}_{s}$.  Hence $V_{s}$ maps
the ideals $\mathfrak{R}_{x}$ and $\mathfrak{R}_{x}^{\star  }$ to
$\mathfrak{R}_{sx}$ and $\mathfrak{R}_{sx}^{\star}$ respectively
since $(\kappa^{\prime}\vert V_{s}\kappa)_{sx} =
(\kappa^{\prime\prime}\vert  \kappa)_{x} = 0$ for
$\kappa\in\mathfrak{R}_{x} $ or $\mathfrak{R}^{\star}   _{x}$.
This gives for each $x$ a well defined linear isometry $V_{s}
\hbox{\rm :}\ \mathfrak{B}/\mathfrak{R}_{x}^{\star}
\mapsto\mathfrak{B}/\mathfrak{R}_{sx}^{\star}$, denoted by
$\mathbf{V}_{s}$, and acting on columns  $\mathbf{k}$ given by the
components
\begin{align*}
(\mathbf{V}_{s}\mathbf{k})_{x}^{-} &= (V_{s}\kappa)_{x}^{-} =
(\delta_{u} \vert V_{s}\kappa)_{sx} =
\sum_{b}l_{sx}(b_{s}^{\star})\kappa^{\ast}_{b} = \kappa_{x}^{-},\\[.3pc]
(\mathbf{V}_{s}\mathbf{k})^{\circ} &= (V_{s}\kappa)^{\circ} = V_{s}   \kappa^{\circ},\\[.4pc]
(\mathbf{V}_{s}\mathbf{k})^{+} &= (V_{s}\kappa)^{+} = \sum_{b
\in\mathfrak{b}_{s}}\kappa_{b^{s}} = \kappa^{+}.
\end{align*}
Hence we may write the action of the semigroup $S$ on
$\mathbb{K}_{x}$ as
\begin{equation*}
\mathbf{V}_{s} =
\begin{bmatrix}  1 & 0 & 0\\[.2pc] 0 & V_{s} & 0\\[.2pc]  0 & 0 & 1
\end{bmatrix},
\end{equation*}
and extend it to $\mathbb{K}_{\Delta}$ as the direct integral.
Clearly it  leaves $\mathbf{e}$ invariant. Finally, we have for
any $\kappa\in\mathfrak{B}$,
\begin{equation*}
(\delta_{b_{s}} \vert V_{s} \kappa)_{sx} = (\delta_{b}
\vert\kappa)_{x},  \qquad V_{s} \delta_{b^{\star}} =
\delta_{b_{s}^{\star}}, \qquad V_{s} (\delta_{b} \cdot\kappa) =
\delta_{b_{s}} \cdot V_{s}\kappa
\end{equation*}
where the first equality follows from the covariance of $l_{x}$,
the second simply from the definition of $V_{s}$, and the third
from the surjectivity of $b \mapsto b^{s}$. This shows that
$\mathbf{V}_{s}\mathbf{j}_{x}
(b)=\mathbf{j}_{sx}(b_{s})\mathbf{V}_{s}$ and hence the required
relation for  integrable subsets $\Delta$.


\begin{thebibliography}{9}
\bibitem{Belavkin85}
Belavkin V P, A reconstruction theorem for a quantum stochastic
process, {\it Teoret. Mat. Fiz.} \textbf{62(3)} (1985) 408--431

\bibitem{Belavkin92c}
Belavkin V P, Chaotic states and stochastic integration in quantum
systems, {\it Russian Mathematical Surveys} \textbf{47(1)} (1992)
47--106

\bibitem{Belavkin03}
Belavkin V P, Quantum {I}t\^{o} algebras: Axioms, representations,
decompositions, {\it Quantum Probability Communications}
\textbf{XI} (2003) 39--54

\bibitem{Bhat94}
Bhat B~V Rajarama and Parthasarathy K~P, Kolmogorov's existence
theorem for {M}arkov processes in ${C}^{\ast}$-algebras, {\it
Proc. Indian Acad.  Sci. (Math. Sci.)} \textbf{104(1)} (1994)
253--262

\bibitem{Chakraborty01}
Chakraborty P S, Goswami D and Sinha K B, A covariant quantum
stochastic dilation theory, {Stochastics in Finite and Infinite
Dimensions} (2001) pp.~89--99

\bibitem{Franz05a}
Franz U and Sch\"urmann M (eds), Quantum independent increment
processes {II}, Lecture Notes in Mathematics, vol.~1866
(Springer-Verlag) (2005)

\bibitem{Heyer95}
Heyer H, Generation and representation of convolution hemigroups
on a  {P}olish group, Infinite-dimensional harmonic analysis
(Tubingen) (1995)

\bibitem{Schurmann93}
Sch\"{u}rmann M, White noise on bialgebras, Lecture Notes in
Mathematics, vol.~1544 (Springer-Verlag) (1993)
\end{thebibliography}
\end{document}